\documentclass[12pt]{iopart}

\usepackage{iopams}

\usepackage[nodots]{numcompress}
\usepackage{verbatim}
\usepackage{graphicx}
\usepackage{fancybox}
\usepackage{subfigure}
\usepackage[usenames,dvipsnames]{color}
\usepackage{sidecap}
\usepackage{caption}
\usepackage{ifthen}
\usepackage{textcomp}
\usepackage[all]{xy}

\newcommand{\N}{{\mathbb N}}

\newcommand{\R}{{\mathbb R}}

\newcommand{\dem}{{\em Proof: \;}}
\newcommand{\fdem}{\hfill $\square$}

\newtheorem{teo}{Theorem}[section]
\newtheorem{lema}[teo]{lemma}
\newtheorem{cor}[teo]{Corollary}

\newtheorem{defi}{Definition}[section]

\begin{document}

\title[Factorization of Constrained Energy Network Reliability]{Factorization of Constrained Energy K-Network Reliability with Perfect Nodes}

\author{Juan Manuel Burgos}

\address{Instituto de Matem\'aticas, Universidad Nacional Aut\'onoma de M\'exico, Unidad Cuernavaca.\\ Av. Universidad s/n, Col. Lomas de
Chamilpa. Cuernavaca, Morelos M\'exico, 62209.}
\ead{burgos@matcuer.unam.mx}

\begin{abstract}
This paper proves a new general $K$-network constrained energy reliability global factorization theorem. As in the unconstrained case, beside its theoretical mathematical importance the theorem shows how to do parallel processing in exact network constrained energy reliability calculations in order to reduce the processing time of this NP-hard problem. Followed by a new simple factorization formula for its calculation, we propose a new definition of constrained energy network reliability motivated by the factorization theorem and the accomplishment of parallel processing, something impossible with the original definition.
\end{abstract}


\maketitle

\section{Introduction}

The energy constraint is a very natural one in real networks what makes the subject of constrained energy $K$-reliability a necessary one in engineering. It is the probability that a state be a $K$-PathSet (or $K$-operative state) with a number of operative edges less than or equal to a given bound. The term energy comes form the prototypical example of a network requiring an amount of energy for each operational edge. We can think of a network whose edges are power lines requiring a cooler device because of the increasing temperature.

We model this situation in the paper, assuming that every edge requires the same amount of energy per unit time to be operational. We normalize to one the energy per unit time consumed by each operational state. In real situations, the amount of energy per unit time consumed by an operational edge, depends on the edge. In the above example, the energy per unit time depends on the edge length (it is proportional neglecting non linear effects). By introducing fictitious nodes to the network, we can approximate the real situation to the one modeled here and the approximation can be made as good as we want just introducing enough nodes (This is described carefully in the author thesis [Bu], chapter 3).

A previous paper by the author [BR] gives a new "global" factorization theorem which allows, in particular, to do parallel processing in order to reduce the computational time in the calculation of this $NP$-hard problem. It is hard to find, if there is any, global factorization graph theorems besides the one in [BR]. Generalizing the graph invariant a little bit makes the problem too hard and in fact the author conjectures that there is no global factorization in most of these cases. Such is the case of the Tutte polynomial and its particular cases or constrained diameter reliability. In this spirit, is very interesting that a new definition of constrained energy reliability makes possible a global factorization of this non trivial and very practical case, factorization which is \textbf{impossible} with the naive definition. The accomplishment of parallel processing the constrained energy reliability is enough justification for introducing this new definition. This new definition is followed by a simple (edge by edge) factorization formula, the analogue of the well known simple factorization in usual reliability [Mo], resulting in a recursive algorithm for its exact calculation. Of course, its approximate calculation can be made by Monte-Carlo methods.

As in the previous paper [BR], the new constrained energy $K$-network reliability factorization theorem gives as particular cases the constrained energy version of the well known reduction transformations (series-parallel, polygon-to-chain [Wo] and delta-star [Ga]) which are the key stone of the known factoring algorithms [SC] for the network reliability exact calculation. As it is mentioned in [BR], besides the well known factorization through an articulation point, no other general "global" factorization theorem is known in exact $K$-network reliability calculation and even less in the constrained energy case. This paper gives a new general "global" factorization theorem solving the following problem:

\textbf{Problem:} \emph{Given a decomposition of a stochastic graph $G$ by subgraphs $G_{1}$ and $G_{2}$ only sharing nodes, express the constrained energy reliability of $G$ in terms of the constrained energy reliabilities of the graphs resulting from $G_{1}$ and $G_{2}$ identifying the common nodes shared by them in all possible ways.}

\section{Preliminaries}

The mathematical model of a Network whose nodes are perfect and its edges can fail is a stochastic graph [Co]; i.e. an undirected graph with associated Bernoulli variables to its edges.

\begin{defi}
An undirected graph $G$ is a pair $(V,E)$ such that $V$ is finite set whose elements will be called nodes and $E$ is a subset of $\{\{a,b\}\ /\ a,b\in V\}\times \N$ whose elements will be called edges such that for each pair of distinct edges $(\{a_{1},b_{1}\},n_{1}), (\{a_{2},b_{2}\},n_{2}) \in E$ we have that $n_{1}\neq n_{2}$.
\end{defi}

\begin{defi}
A stochastic graph $G$ is a tern $(V,E, \Phi)$ such that $(V, E)$ is a graph and $\Phi:E\rightarrow Ber$ is a function which associates a Bernoulli variable to each edge in such a way that these variables are independent.
\end{defi}

Each Bernoulli variable is characterized by a parameter $p$ in the $[0,1]$ closed interval and we can write a stochastic graph as $(G, \{p_{e}\}_{e\in E})$ where $G$ is an undirected graph and $p_{e}$ is the parameter of the variable $\Phi(e)$. Nodes and edges of $G$ will be denoted by $V(G)$ and $E(G)$ respectively.

\begin{defi}
A state $\mathcal{E}$ of the graph $G=(V,E)$ is a function $\mathcal{E}: E\rightarrow \{0, 1\}$. An edge $e$ will be called operative if $\mathcal{E}(e)=1$ and will be called non-operative otherwise.
\end{defi}

Consider a subset $K$ of $V(G)$. A state $\mathcal{E}$ of the graph $G$ will be called a $K$-PathSet (or $K$-operative) if $K$ is contained in the set of nodes of any of the edge-connected components of the graph resulting from removing the non-operative edges of $G$. Otherwise the state will be called a $K$-CutSet. We will denote by $\#\mathcal{E}$ the number of operative edges of $\mathcal{E}$; i.e. $$ \# \mathcal{E}= \# \mathcal{E}^{-1}(1)$$

\begin{defi}\label{DefConfiabilidad}
The constrained $l$-energy $K$-reliability of a stochastic graph $G$ is $$R_{K,l}(G)=P(\mathcal{E}\ is\ a\ K-PathSet\ and\ \#\mathcal{E}\leq l)$$
\end{defi}

Because of the independence of the Bernoulli variables associated to the edges, we can calculate the constrained $l$-energy $K$-reliability in the following way:

\begin{equation}\label{ProbabilidadPathset}
    P(\mathcal{E})=\prod_{e_{i}\in E(G)}p_{i}^{\mathcal{E}(e_{i})}(1-p_{i})^{1-\mathcal{E}(e_{i})}
\end{equation}

$$R_{K,l}(G)=\sum_{\mathcal{E}\ is\ a\ K-PathSet\ and\ \#\mathcal{E}\leq l} P(\mathcal{E})$$

The resulting algorithm from the above expression is uneffective for it requires a complete list of the operative states involved. In the next section a much more effective recursion algorithm for the exact calculation of the constrained energy reliability is shown.

\section{Constrained Energy Reliability}

\begin{figure}
\begin{center}
  \includegraphics[width=0.4\textwidth]{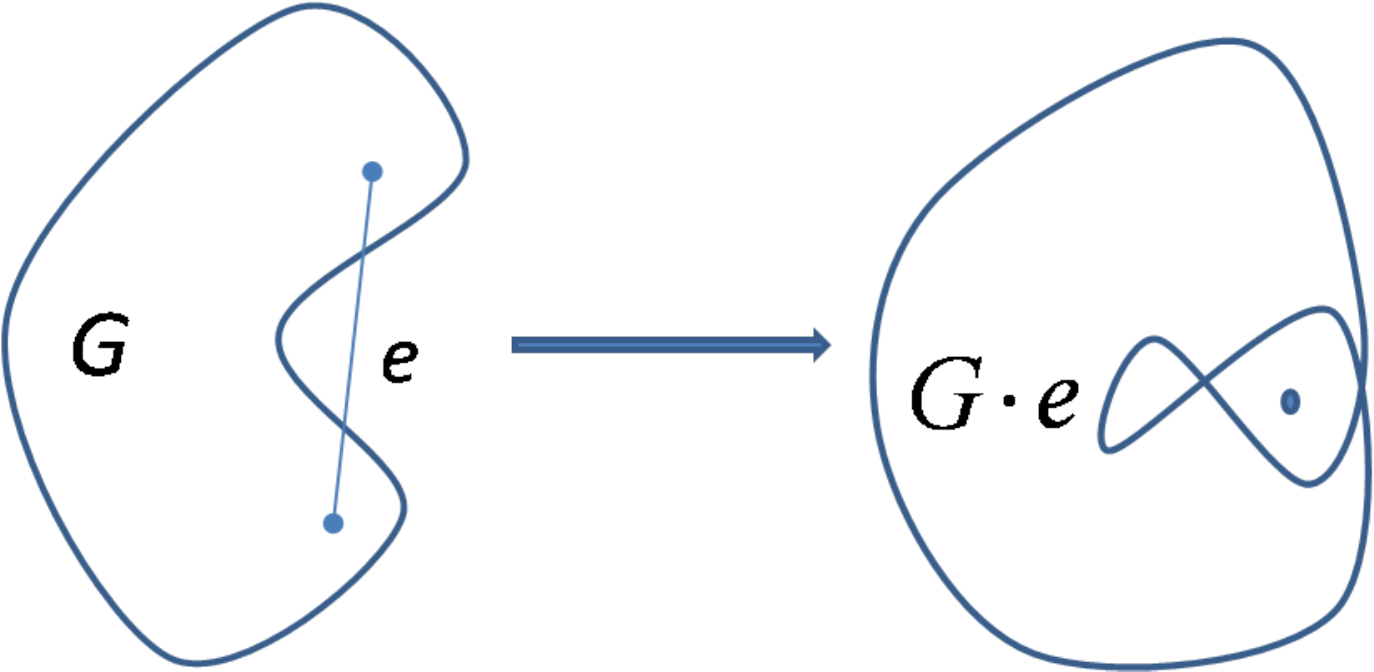}\\
  \end{center}
  \caption{Contraction of the edge \textit{e}}\label{contrac}
\end{figure}

\begin{defi}
Consider an edge $e=(\{v,w\}, n) \in E$ of a graph $G=(V,E)$ and define the following equivalence relation in $V$: $a\sim b$ if $a=b$ or $\{a,b\}=\{v,w\}$. Consider the suryective canonical function $\pi:V\rightarrow V/\sim$ such that $\pi(a)=[a]_{\sim}$. We define the contraction of an edge $e$ in $G$ as the graph $G\cdot e$ such that $$G\cdot e = (V/\sim, E\cdot e)$$ where $E\cdot e =\{\ (\{\pi(a),\pi(b)\},n)\ /\ (\{a,b\},n)\in E-\{e\}\ \}$ (see Figure \ref{contrac}). We will denote by $K_{e}=\pi(K)$ the new distinguished set of nodes contained in $V(G\cdot e)$ where $K$ is the distinguished set of nodes in $V(G)$.
\end{defi}

\begin{figure}
\begin{center}
  \includegraphics[width=0.4\textwidth]{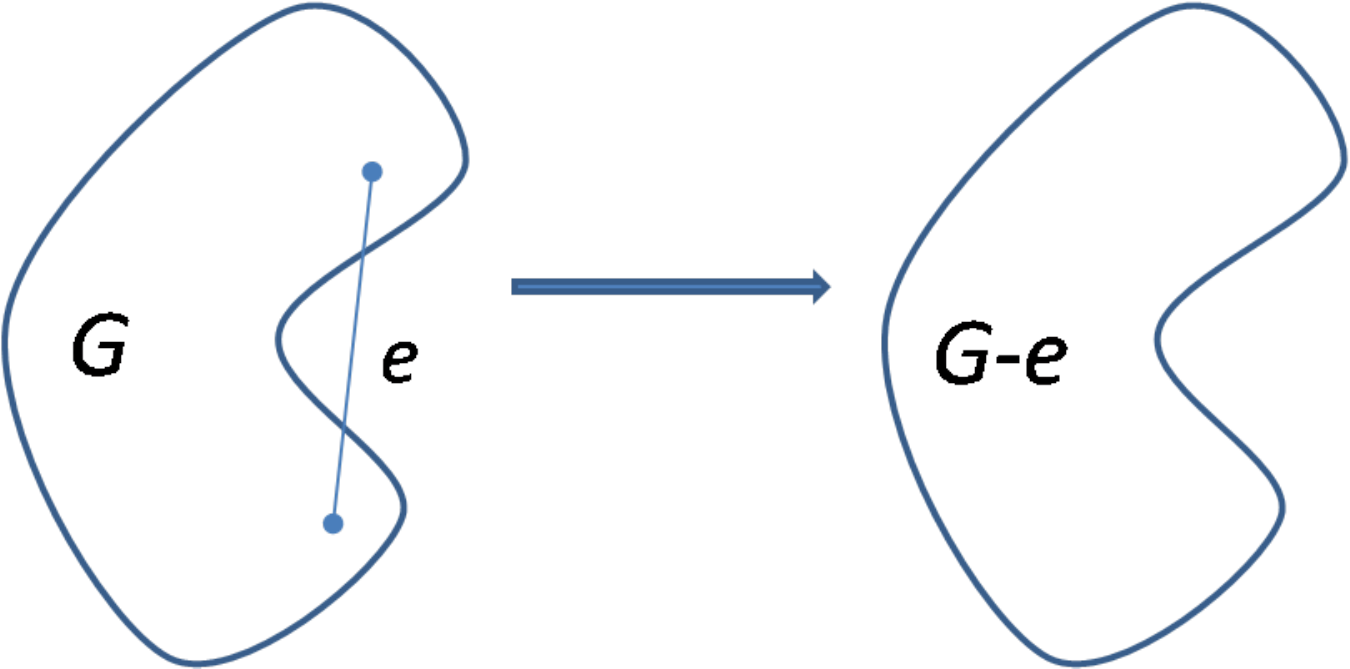}\\
  \end{center}
  \caption{Deletion of the edge \textit{e}}\label{extrac}
\end{figure}

\begin{defi}
Consider an edge $e=(\{v,w\}, n) \in E$ of the graph $G=(V,E)$. We define the deletion of the edge $e$ of $G$ as the graph (see Figure \ref{extrac}) $$G-e=(V,E-\{e\}\ )$$
\end{defi}

\begin{defi}
Define the following partial order in the set of states of $G$: $\mathcal{E} \leq \mathcal{F}$ if $\mathcal{E}^{-1}(1) \subset \mathcal{F}^{-1}(1)$. The state $\mathcal{E}$ is a $K$-minpath if it is minimal in the set of $K$-PathSets.
\end{defi}

Is clear that if the given energy bound is strictly less than the number of operative edges in each $K$-minpath, then the constrained energy reliability is zero; i.e There is no enough energy to turn on the network.

\begin{defi}
Given a stochastic graph $G$ we define its $\mathcal{I}(G)\in \R[x]$ polynomial as $$\mathcal{I}(G)= \sum_{i=0}^{+\infty}p_{i}x^{i}$$ such that
$$a_{i}=\sum_{\mathcal{E}\ is\ a\ K-PathSet\ and\ \#\mathcal{E}= i} P(\mathcal{E})$$
\end{defi}
Is clear that $$\mathcal{I}(G)= \sum_{i=m}^{n}a_{i}x^{i}$$ such that $n$ is the number of edges of $G$ and $m$ is the greatest lower bound of the number of edges of $K$-PathSets; i.e. The threshold energy of the network.

\begin{figure}
\begin{center}
  \includegraphics[width=0.8\textwidth]{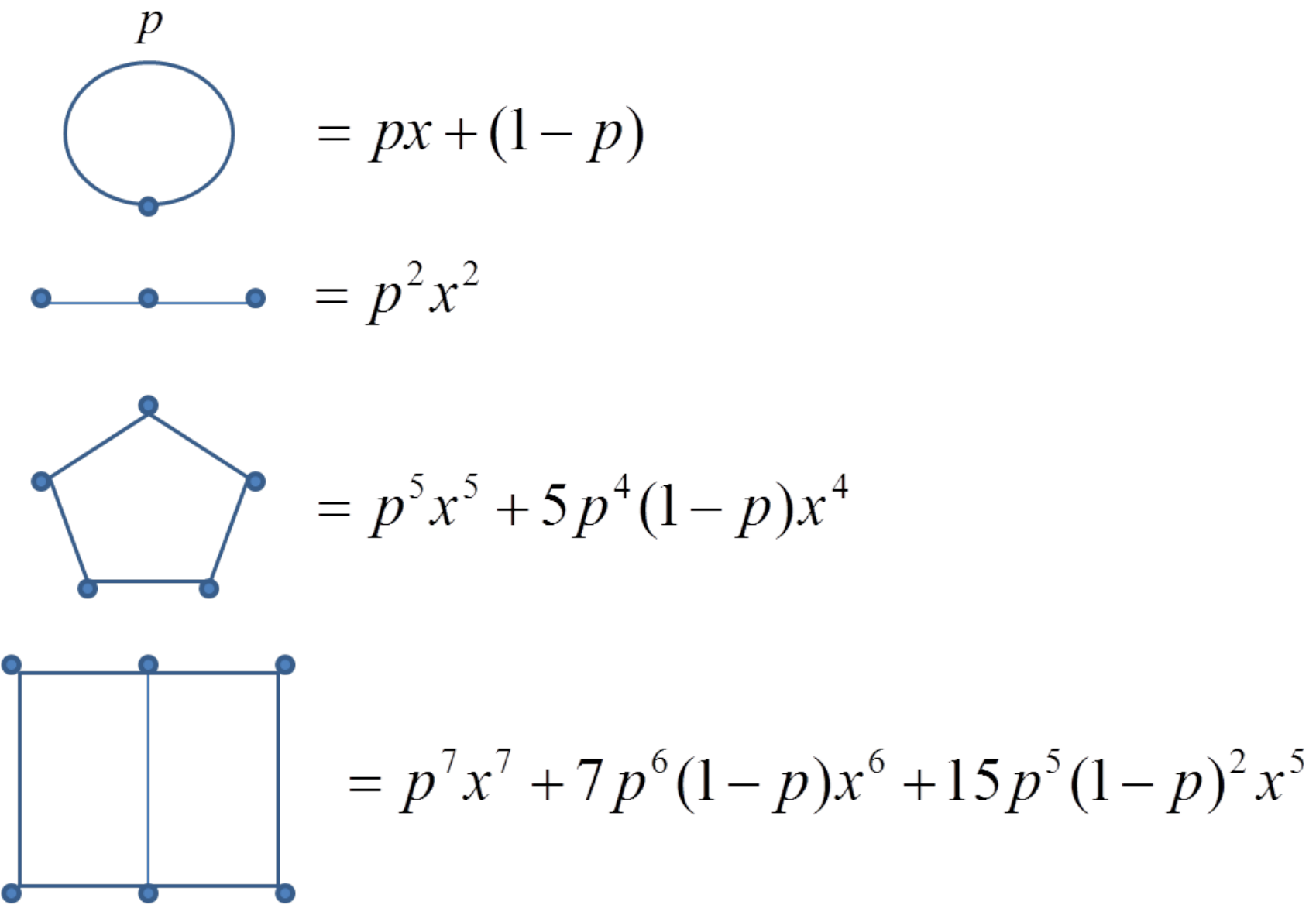}\\
  \end{center}
  \caption{Examples of I(G)}\label{Examples}
\end{figure}

In Figure \ref{Examples} some examples of $\mathcal{I}(G)$ for homogeneous graphs (equal edge reliability $p$) are given. In particular, from the first example from above, we see that irrelevant edges for the calculation of the usual reliability $R_{K}(G)$ are \textbf{relevant} in the calculation of $\mathcal{I}_{K}(G)$. This makes sense because even an irrelevant edge (in terms of $K$-reliability) consumes energy if it is operational and this fact must be considered in the formalism.

The following is the simple factorization formula for the calculation of $\mathcal{I}_{K}(G)$.

\begin{lema}\label{SimpleFact_I}
Consider a stochastic graph $G$ and a subset $K$ of its nodes. For any edge $e_{i}$ of $G$ we have that $$\mathcal{I}_{K}(G)=p_{i}x\ \mathcal{I}_{K_{i}}(G\cdot e_{i}) + (1-p_{i})\mathcal{I}_{K}(G-e_{i})$$
\end{lema}
\dem
Considering arbitrary edge probabilities $(p_{e})_{e\in E}$ of the stochastic graph $G=(V,E,(p_{e})_{e\in E})$, we can see the $K$-reliability $R_{K}(G)$ as a polynomial in the \textbf{formal} variables $(p_{e})_{e\in E}$ and $(q_{e})_{e\in E}$; i.e. $$R_{K}(G)\in \R[(p_{e})_{e\in E}, (q_{e})_{e\in E}]$$
$$R_{K}(G)=\sum_{\mathcal{E}\ is\ a\ K-PathSet}\left( \prod_{e_{i}\in E(G)}p_{i}^{\mathcal{E}(e_{i})}q_{i}^{1-\mathcal{E}(e_{i})}\right)$$
Almost verbatim we can adapt the proof of the well known [Mo] simple factorization in the $(p_{e})_{e\in E}$ variables for the simple factorization in the $(p_{e})_{e\in E}$ and $(q_{e})_{e\in E}$ variables:
$$R_{K}(G)=p_{i}R_{K_{i}}(G\cdot e_{i}) + q_{i}R_{K}(G-e_{i})$$
Because the variables in the polynomial are \textbf{formal}, we can just replace the variable $p_{e}$ for $p_{e}x$ where $x$ is a new formal variable. Is clear that the resulting polynomial $$\hat{\mathcal{I}}_{K}(G)\in \R[x, (p_{e})_{e\in E}, (q_{e})_{e\in E}]$$ is $$\hat{\mathcal{I}}_{K}(G)=\sum_{i=m}^{n} a_{i}x^{i}$$ such that
$$a_{i}=\sum_{\mathcal{E}\ is\ a\ K-PathSet\ and\ \#\mathcal{E}= i} \left(\prod_{e_{i}\in E(G)}p_{i}^{\mathcal{E}(e_{i})}q_{i}^{1-\mathcal{E}(e_{i})}\right)$$
and verifies $$\hat{\mathcal{I}}_{K}(G)=p_{i}x\ \hat{\mathcal{I}}_{K_{i}}(G\cdot e_{i}) + q_{i}\hat{\mathcal{I}}_{K}(G-e_{i})$$
Evaluating the variables $(q_{e})_{e\in E}$ by $q_{e}= 1-p_{e}$ and denoting the resulting polynomial by $\mathcal{I}_{K}(G)$, we have the lemma.
\fdem

Is clear from the definition that the $K$-reliability of $G$ is just the evaluation of $\mathcal{I}(G)\in \R[x]$ in $x=1$: $$R_{K}= ev_{1}\circ \mathcal{I}_{K}$$ This motivates the consideration of the $l$-truncated polynomial $$\mathcal{I}_{K,l}= \sum_{i=m}^{l} a_{i}x^{i}$$ such that
$$a_{i}=\sum_{\mathcal{E}\ is\ a\ K-PathSet\ and\ \#\mathcal{E}= i} P(\mathcal{E})$$
This way we have that $$R_{K,l}= ev_{1}\circ \mathcal{I}_{K,l}$$
Because of the previous lemma, is clear that the truncated polynomials verify the relation: $$\mathcal{I}_{K,l}(G)=p_{i}x\ \mathcal{I}_{K_{i},l-1}(G\cdot e_{i}) + (1-p_{i})\mathcal{I}_{K,l}(G-e_{i})$$
We have proved the following simple factorization for energy constrained reliability:

\begin{cor}
Consider a stochastic graph $G$ and a subset $K$ of its nodes. For any edge $e_{i}$ of $G$ we have that $$R_{K,l}(G)=p_{i}R_{K_{i},l-1}(G\cdot e_{i}) + (1-p_{i})R_{K,l}(G-e_{i})$$
\end{cor}

Because of the factorization theorem we are looking for, we propose the following definition for constrained energy reliability:

\begin{defi}
Consider a stochastic graph $G$ and a subset $K$ of its nodes. The constrained $l$-energy $K$-reliability of $G$ is $$[\mathcal{I}_{K}(G)]_{l+1}\in \frac{\R[x]}{\langle x^{l+1}\rangle}$$
\end{defi}

Is clear that the $l$-truncated polynomial $\mathcal{I}_{K,l}(G)$ is a representative of $[\mathcal{I}_{K}(G)]_{l+1}$; i.e
$$\mathcal{I}_{K,l}(G) + \langle x^{l+1}\rangle = [\mathcal{I}_{K}(G)]_{l+1}$$
Moreover, $\mathcal{I}_{K,l}(G)$ is the only representative from which we can calculate the original constrained energy reliability. This is the relation between the original definition and the one proposed here. From here to the rest of the paper, constrained energy reliability means our new definition.

Because $\mathcal{I}_{K,l}(G)$ is a representative of $[\mathcal{I}_{K}(G)]_{l+1}$, the calculation of the exact constrained energy reliability can be performed via the resulting recursive algorithm of the expression $$\mathcal{I}_{K,l}(G)=p_{i}x\ \mathcal{I}_{K_{i},l-1}(G\cdot e_{i}) + (1-p_{i})\mathcal{I}_{K,l}(G-e_{i})$$ As an approximate method, we can approximate $\mathcal{I}_{K}(G)$ via a Monte-Carlo method and then truncate the approximate polynomial at the given energy bound.

\section{The Factorization Theorem}

\begin{flushleft}
\textbf{Hypothesis 1:} In the whole paper, $G_{1}$, $G_{2}$ and $G$ are graphs with distinguished subset of nodes $K_{1}$, $K_{2}$ and $K$ respectively such that $G=G_{1}\cup G_{2}$, $K=K_{1}\cup K_{2}$ and
\end{flushleft}
$$\{k_{1},k_{2},\ldots k_{n}\}=K_{1}\cap K_{2}=G_{1}\cap G_{2}$$
\textbf{Hypothesis 2:} We assume in the whole paper that for each node $v\in K$ there exists a path in $G$ that joins $v$ with some vertex $k_{i}\in K_{1}\cap K_{2}$.\\

\begin{figure}
\begin{center}
  \includegraphics[width=0.3\textwidth]{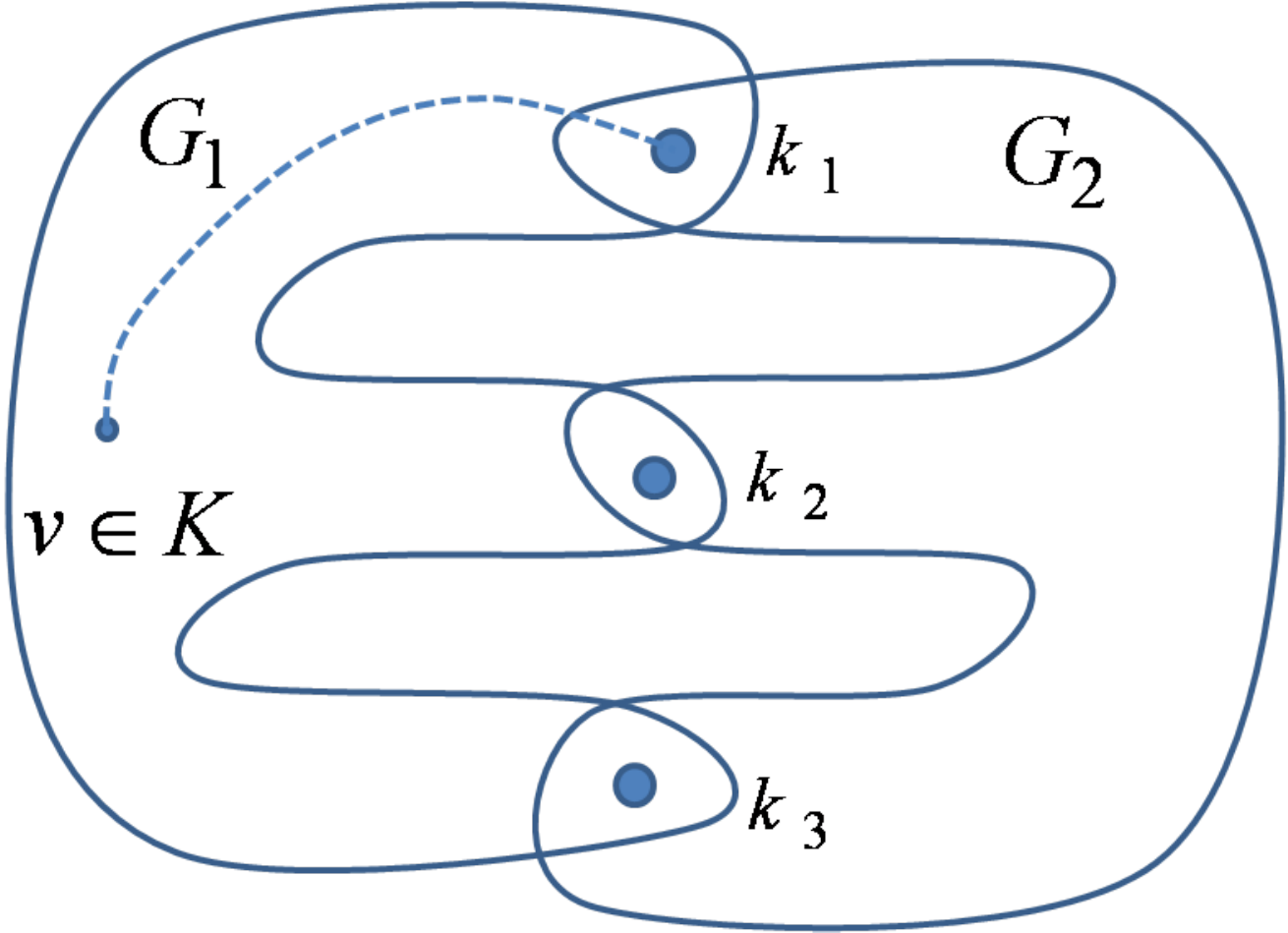}\\
\end{center}
\caption{Hypothesis on the Graph}\label{Hypothesis}
\end{figure}

These Hypothesis are illustrated in figure \ref{Hypothesis}. In view of the first hypothesis, it is reasonable to assume the second one, otherwise $R_{K}(G)=0$ and there would be no necessity for any calculation. For notational convenience, the $K$ subscript in $R_{K}(G)$ will be omitted in the rest of the paper.

\begin{lema}
$G$ is $K$-connected if and only if $G$ is $\{k_{1},k_{2},\ldots k_{n}\}$-connected.
\end{lema}
\dem
The direct is trivial. Conversely, take a pair o nodes $a$ and $b$ in $K$. There are paths $P_{a}$ y $P_{b}$ connecting $a$ and $b$ with $k_{i}$ y $k_{j}$ respectively. Because $G$ is $\{k_{1},k_{2},\ldots k_{n}\}$-connected, there is a path $P$ connecting $k_{i}$ with $k_{j}$. The concatenation of the paths $P_{a}$, $P$ and $P_{b}$ joins $a$ with $b$.
\fdem

The previous lemma motivates the following definition.

\begin{defi}
Consider the equivalence relation: $k_{i}\sim^{l} k_{j}$ if there is a path in  $G_{l}$ joining $k_{i}$ with $k_{j}$, $l=1,2$.
The connectivity state of $G_{l}$ is the partition of $\{k_{1},k_{2},\ldots k_{n}\}$ given by $$\mathcal{C}_{l}=\{k_{1},k_{2},\ldots k_{n}\}/\sim^{l}$$
\end{defi}

\begin{figure}
\begin{center}
  \includegraphics[width=0.6\textwidth]{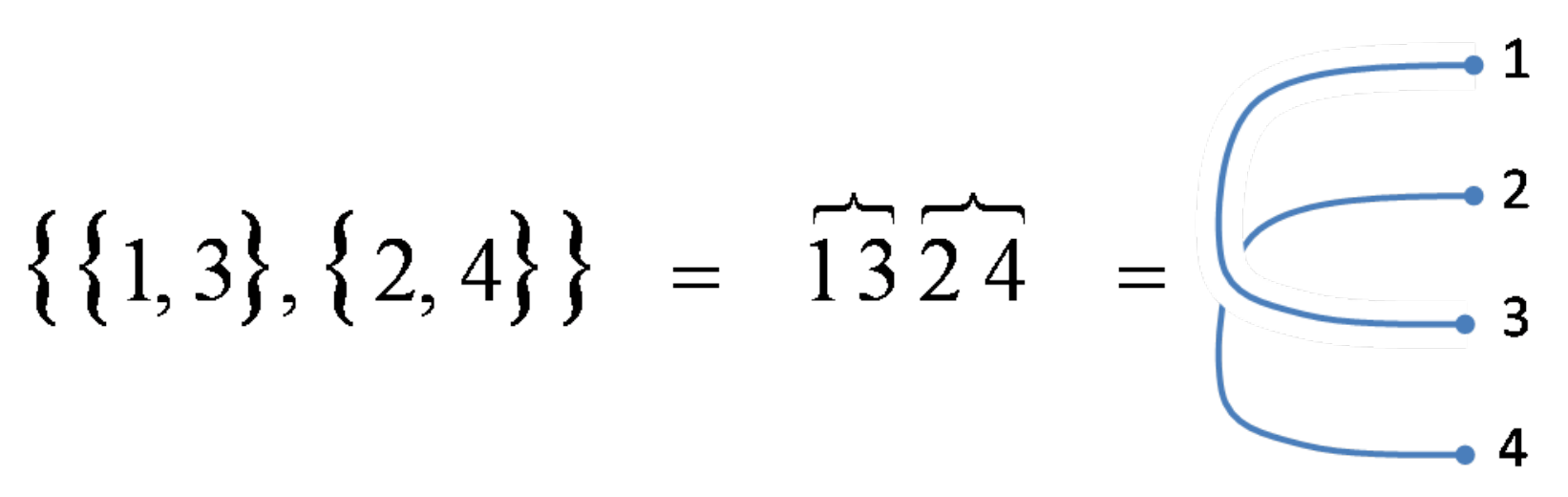}\\
  \end{center}
  \caption{Connectivity State}\label{EjEstadoConect}
\end{figure}

Denote by $Con$ the set of partitions of $\{k_{1},k_{2},\ldots k_{n}\}$ whose elements will be called connectivity states. Figure \ref{EjEstadoConect} shows some useful notational and diagrammatical  ways to represent a connectivity state.

\begin{defi}
For each connectivity state $\mathcal{A}$ denote by $G^{\mathcal{A}}_{l}$ the graph resulting of the identification of the nodes in $\{k_{1},k_{2},\ldots k_{n}\}$ of $G_{l}$ by the state $\mathcal{A}$; i.e. given the graph $G_{l}=(V,E)$ define $G^{\mathcal{A}}_{l} = (V/\sim^{\mathcal{A}}_{l}, E^{\mathcal{A}})$ where $$E^{\mathcal{A}}= \{\ (\{\pi(a),\pi(b)\},n)\ /\ (\{a,b\},n)\in E\ \}$$ and $\sim^{\mathcal{A}}_{l}$ is the equivalence relation in $V$ generated by $\mathcal{A}$ with the canonical suryection $$\pi:V\rightarrow V/\sim^{\mathcal{A}}_{l}$$ such that $\pi(a)=[a]_{\sim^{\mathcal{A}}_{l}}$. The distinguished set of nodes is $K^{\mathcal{A}}_{l}= \pi(K_{l})$.
\end{defi}

\begin{defi}
The set of connectivity states of $G$ is $Con\times Con$ and $(\mathcal{C}_{1}, \mathcal{C}_{2})$ is the connectivity state of $G$ where $\mathcal{C}_{l}$ is the connectivity state of $G_{l}$, $l=1,2$.
\end{defi}

\begin{defi}
We say that a connectivity state $(\mathcal{A}, \mathcal{B})$ of $G$ is connected if $$\{k_{1},k_{2},\ldots k_{n}\}/\sim = \{\{k_{1},k_{2},\ldots k_{n}\}\}$$ where $\sim$ is the following equivalence relation in $\{k_{1},k_{2},\ldots k_{n}\}$: $k_{i}\sim k_{j}$ if $k_{i}\sim_{\mathcal{A}} k_{j}$ or $k_{i}\sim_{\mathcal{B}} k_{j}$.
\end{defi}

\begin{figure}
\begin{center}
  \includegraphics[width=0.4\textwidth]{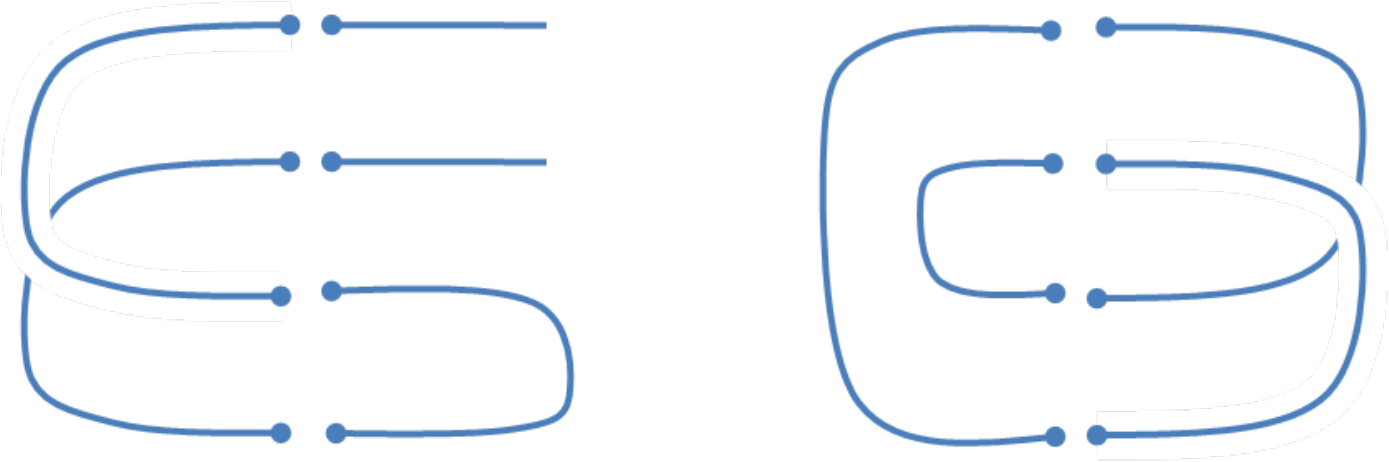}\\
  \end{center}
  \caption{Connected connectivity states of $G$}\label{EjEstadoConexo}
\end{figure}

\begin{figure}
\begin{center}
  \includegraphics[width=0.4\textwidth]{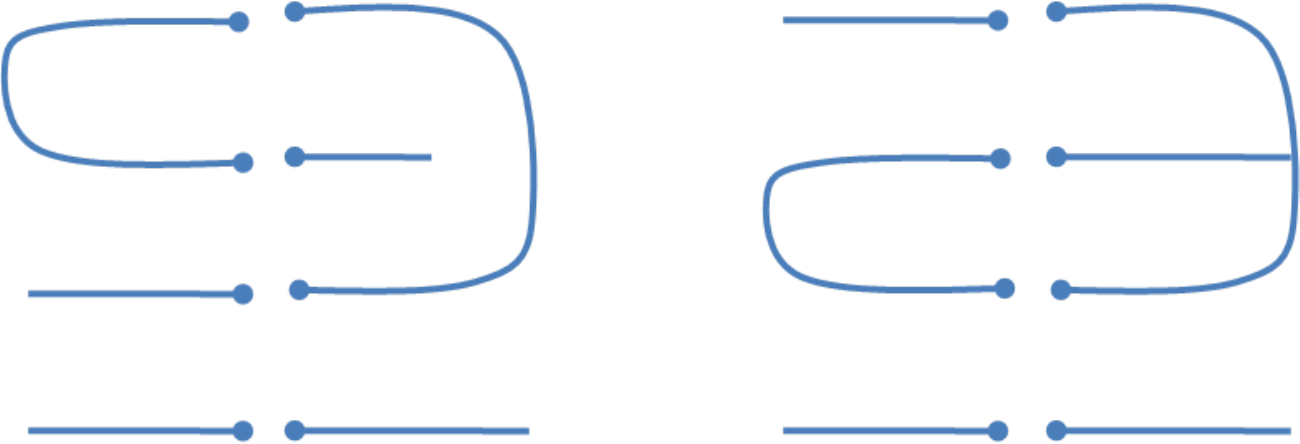}\\
  \end{center}
  \caption{Non connected connectivity states of $G$}\label{EjEstadoNoConexo}
\end{figure}

Figures \ref{EjEstadoConexo} and \ref{EjEstadoNoConexo} show examples of connectivity states of $G$.

Considering an ordering in $Con$ we define the connectivity matrix $A=(a_{ij})$ given by $a_{ij}=1$ if $(\mathcal{A}_{i}, \mathcal{A}_{j})$ is connected and $a_{ij}=0$ if it is not. For example, in the case of three sharing nodes, ordering the base $Con$ in the following way: $$Con= \{123, 1\overbrace{23}, \overbrace{13}2, \overbrace{12}3, \overbrace{123}\}$$ we get the connectivity matrix $$A=\left(
                                                                \begin{array}{ccccc}
                                                                  0 & 0 & 0 & 0 & 1 \\
                                                                  0 & 0 & 1 & 1 & 1 \\
                                                                  0 & 1 & 0 & 1 & 1 \\
                                                                  0 & 1 & 1 & 0 & 1 \\
                                                                  1 & 1 & 1 & 1 & 1 \\
                                                                \end{array}\right)$$
One of the deep results in [BR] is the fact that the connectivity matrix is invertible. Moreover, a beautiful formula for its determinant is given.

\begin{lema}
Let $(b_{ij})=A^{-1}$ where $A$ is the connectivity matrix. Then $$\mathcal{I}_{K}(G)=\sum_{i,j=1}^{m}b_{ij}\ \mathcal{I}_{K_{i}}(G_{1}^{\mathcal{A}_{i}}) \ \mathcal{I}_{K_{j}}(G_{2}^{\mathcal{A}_{j}})$$ and the above expression doesn't depend on the order of the base $Con$.
\end{lema}
\dem
In the same fashion as in lemma \ref{SimpleFact_I}, we can see $K$-reliability as a polynomial in the \textbf{formal} variables $(p_{e})_{e\in E}$ and $(q_{e})_{e\in E}$; i.e. $$R_{K}(G)\in \R[(p_{e})_{e\in E}, (q_{e})_{e\in E}]$$ Almost verbatim, we can prove the factorization theorem with the proof given in [BR] this time in the euclidean domain $\R[(p_{e})_{e\in E}, (q_{e})_{e\in E}]$ instead in the field $\R$ (we just have to consider now the functional 
$$P_{1}\otimes P_{2}:\ V\otimes V \rightarrow \R[(p_{e})_{e\in E}, (q_{e})_{e\in E}]$$ which translates the combinatorial problem into the probabilistic problem in [BR]):
$$R_{K}(G)=\sum_{i,j=1}^{m}b_{ij}\ R_{K_{i}}(G_{1}^{\mathcal{A}_{i}})\ R_{K_{j}}(G_{2}^{\mathcal{A}_{j}})$$ and the above expression doesn't depend on the order of the base $Con$. Because the variables in the polynomial are \textbf{formal}, we can just replace the variable $p_{e}$ for $p_{e}x$ where $x$ is a new formal variable. Following the notation in the proof of lemma \ref{SimpleFact_I} we have in the euclidean domain $\R[x, (p_{e})_{e\in E}, (q_{e})_{e\in E}]$ the following identity: $$\hat{\mathcal{I}}_{K}(G)=\sum_{i,j=1}^{m}b_{ij}\ \hat{\mathcal{I}}_{K_{i}}(G_{1}^{\mathcal{A}_{i}}) \ \hat{\mathcal{I}}_{K_{j}}(G_{2}^{\mathcal{A}_{j}})$$ and the above expression doesn't depend on the order of the base $Con$. We have shown in the mentioned lemma that evaluating the variables $(q_{e})_{e\in E}$ by $q_{e}= 1-p_{e}$ in $\hat{\mathcal{I}}_{K}(G)$ gives the polynomial $\mathcal{I}_{K}(G)$. Because the evaluation is an algebra morphism, we have the result.
\fdem

Because $\langle x^{l+1}\rangle \lhd \R[x]$; i.e. $\langle x^{l+1}\rangle$ is an ideal of the polynomial ring, the algebra structure in the quotient can be defined such that the canonical epimorphism $\pi:\R[x]\rightarrow \frac{\R[x]}{\langle x^{l+1}\rangle}$ is an algebra morphism; i.e. $$[a]_{l+1}[b]_{l+1}=[ab]_{l+1}$$ where $a$ and $b$ are polynomials in $\R[x]$. We have proved the main theorem of the paper:

\begin{teo}
Let $(b_{ij})=A^{-1}$ where $A$ is the connectivity matrix. Then $$[\mathcal{I}_{K}(G)]_{l+1}=\sum_{i,j=1}^{m}b_{ij}\ [\mathcal{I}_{K_{i}}(G_{1}^{\mathcal{A}_{i}})]_{l+1} \ [\mathcal{I}_{K_{j}}(G_{2}^{\mathcal{A}_{j}})]_{l+1}$$ and the above expression doesn't depend on the order of the base $Con$.
\end{teo}

\begin{figure}
\begin{center}
  \includegraphics[width=1.0\textwidth]{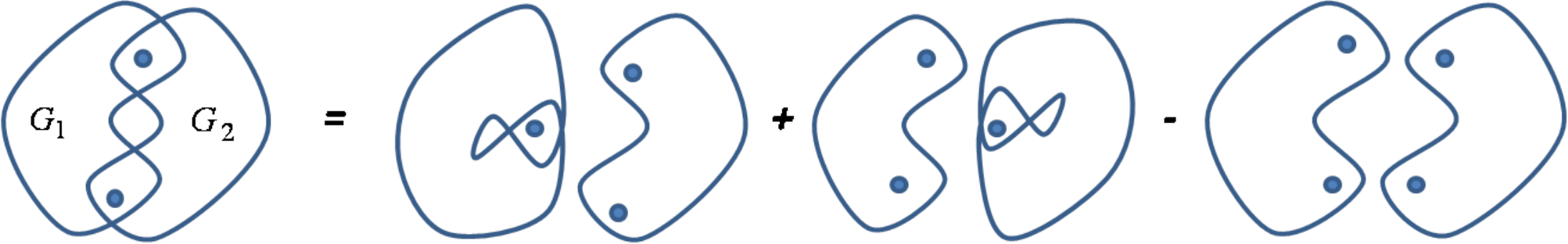}\\
  \end{center}
  \caption{}\label{nDos}
\end{figure}

\begin{figure}
\begin{center}
  \includegraphics[width=1.0\textwidth]{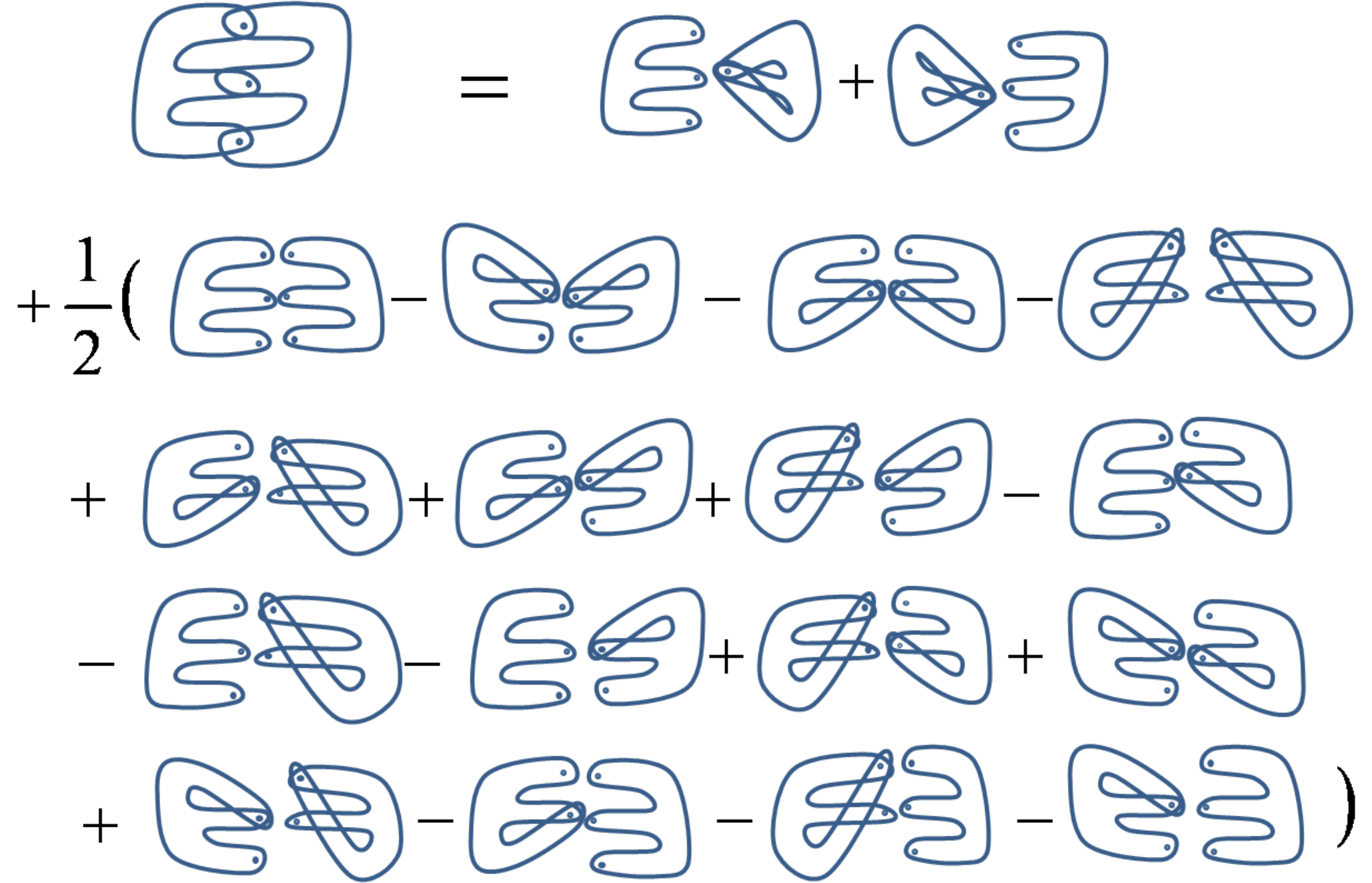}\\
  \end{center}
  \caption{}\label{nTres}
\end{figure}

Figures \ref{nDos} and \ref{nTres} illustrates the factorization with two and three sharing nodes respectively. Is clear that the theorem is false if we consider the truncated polynomials $\mathcal{I}_{K,l}(G)$ instead of the classes $[\mathcal{I}_{K}(G)]_{l+1}$. In particular, the theorem proves that it is \textbf{impossible} to have a factorization theorem in terms of the original definition of constrained energy reliability. This justifies our definition. Because of the isomorphism: $$\frac{\R[x]}{\langle x^{l+1}\rangle} \simeq \R[\sigma]$$ such that $\sigma^{l+1}=0$ and $\sigma^{l} \neq 0$, to get the actual probability (the original definition) we can do the calculation of the factorization formula with the truncated polynomials $\mathcal{I}_{K,l}(G)$ neglecting every term whose degree is greater than the energy bound and finally evaluate the obtained result in the way explained before. 

For more details on the factorization, the references are [BR] and the author thesis [Bu].

\end{document}